\documentclass{amsart}
\usepackage{graphicx}
\usepackage{amsmath}
\usepackage{amsfonts}
\usepackage{amssymb}
\usepackage{amsthm}
\usepackage{color}
\begin{document}

\def\R{\mathbb{R}}
\def\N{\mathbb{N}}
\def\H{\mathcal{H}}
\def\d{\textrm{div}}
\def\v{\textbf{v}}
\def\I{\hat{I}}
\def\B{\hat{B}}
\def\x{\hat{x}}
\def\y{\hat{y}}
\def\p{\hat{\phi}}
\def\r{\hat{r}}
\def\w{\textbf{w}}
\def\u{\textbf{u}}
\def\K{\mathcal{K}}
\def\Reg{\textrm{Reg}}
\def\s{\textrm{Sing}}
\def\sgn{\textrm{sgn}}

\newtheorem{defin}{Definition}[section]
\newtheorem{lem}{Lemma}[section]
\newtheorem{rem}{Remark}[section]
\newtheorem{cor}{Corollary}[section]
\newtheorem{thm}{Theorem}[section]
\newtheorem{prop}{Proposition}[section]
\newtheorem{definition}{Definition}[section]
\newtheorem{con}{Conjecture}[section]
\newtheorem{Main}{Main Result}

\parskip 3pt

\numberwithin{equation}{section}

\title[Fractional optimal maximization problem]{Fractional optimal maximization problem and the unstable fractional obstacle problem}

\author[J.F. Bonder, Z. Cheng and H. Mikayelyan]{Juli\'an Fern\'andez Bonder,  Zhiwei Cheng and Hayk Mikayelyan}

\address[J.F. Bonder]{Departamento de Matem\'atica FCEN - Universidad de Buenos Aires and IMAS - CONICET. Ciudad Universitaria, Pabell\'on I (C1428EGA) Av. Cantilo 2160. Buenos Aires, Argentina.}

\email{jfbonder@dm.uba.ar}

\urladdr{http://mate.dm.uba.ar/~jfbonder}

\address[Z. Cheng and H. Mikayelyan]{Mathematical Sciences, University of Nottingham Ningbo China, 199 Taikang East Road, Ningbo 315100, PR China.}

\email[H. Mikayelyan]{Hayk.Mikayelyan@nottingham.edu.cn}
\urladdr{https://www.nottingham.edu.cn/en/science-engineering/staffprofile/hayk-mikayelyan-.aspx}

	

\begin{abstract}
We consider an optimal rearrangement maximization problem involving the fractional Laplace operator $(-\Delta)^s$, $0<s<1$, and the Gagliardo-Nirenberg seminorm $[u]_s$. We prove the existence of a maximizer, analyze its properties and show that it satisfies the unstable fractional obstacle problem equation for some $\alpha>0$
$$(-\Delta)^s u=\chi_{\{u>\alpha\}}.$$ 
\end{abstract}


\maketitle

\section{Introduction}

One of the classical problems in rearrangement theory is the maximization of the
functional
\begin{equation*}
\Phi(f)=\int_D |\nabla u_f|^2 dx,
\end{equation*}
where $u_f$ is the unique solution of the Dirichlet boundary value problem
\begin{equation*}
\begin{cases} -\Delta u_f(x) = f (x)  & \mbox{in } D,
 \\
 u_f=0  & \mbox{on } \partial D,  \end{cases}
\end{equation*}
and $f$ belongs to the set
$$
\bar{\mathcal{R}}_\beta=\left\{f\in L^\infty(D)\colon 0\leq f \leq 1,\,\,\int_D fdx=\beta \right\}\subset L^\infty(D),
$$
where $\bar{\mathcal{R}}_\beta$ is the weak-closure of the rearrangement class
$$
\mathcal{R}_\beta=\left\{f\in L^\infty(D)\colon f=\chi_E,\ |E|=\beta \right\}.
$$

The problem and its variations, such as the minimization problem and its $p-$harmonic and constraint cases, has been studied by various authors (see \cite{B1, BM, EL, M, Kbook}), and the results, for this particular setting, can be formulated in the following theorem:
 \begin{thm}
 There exists a solution $\hat{f}\in\mathcal{R}_\beta$ such that
 $$
 \Phi(f)\leq\Phi(\hat{f})
$$
for any $f\in  \bar{\mathcal{R}}_\beta$. Moreover, there exists a constant $\alpha>0$ such that
 $$\hat{f}=\chi_{\{\hat{u}>\alpha\}},$$
 where $\hat{u}=u_{\hat{f}}$.
\end{thm}

Let us observe that as a result the function $U=\alpha- \hat{u}$ will be a solution of the unstable obstacle problem
$$
-\Delta u=\chi_{\{u>0\}},
$$
which is one of the classical free boundary problems (see \cite{MW}).

In this paper we consider the fractional analogue of the optimal rearrangement problem and show that its  
maximizers solve the fractional unstable obstacle problem that was recently consider in \cite{AG}.

For the minimization problem, in \cite{BCM_min} we analyzed the fractional version of the optimal rearrangement minimization and show its connection with the stable fractional free boundary problem. 

Our main result is the following theorem. The reader unfamiliar with the fractional vocabulary can find its basic 
objects, their definitions and properties is Section \ref{sec:fr}. 

Let $0<s<1$ be fixed. To avoid extra notations from now on we will use $u_f$ to denote the solution to 
\begin{equation*}
\begin{cases} (-\Delta)^s u_f(x) = f (x)  & \mbox{in } D,
 \\
 u_f=0  & \mbox{in }  D^c,  \end{cases}
\end{equation*}
and 
$$
\Phi_s(f)=[u_f]_s^2,
$$
where $[u]_s$ is the Gagliardo-Nirenberg semi-norm (see Section \ref{sec:fr}).

The main result of the paper is the following:
\begin{thm}\label{thm:main}
There exists a maximizer $\hat{f}\in \mathcal{R}_\beta$ such that 
$$
 \Phi_s(f)\leq\Phi_s(\hat{f})
$$
for any $f\in  \bar{\mathcal{R}}_\beta$. Moreover, for any maximizer $\hat{f}\in \bar{\mathcal{R}}_\beta$ of $\Phi_s$ 
there exists $\alpha>0$ such that 
$$
\hat{f}=\chi_{\{\hat{u}> \alpha\}},
$$
where $\hat{u}=u_{\hat{f}}$.

As a result the function $\hat{u}$ solves the fractional unstable obstacle equation
\begin{equation}\label{cond-ineq}
(-\Delta)^s \hat{u}= \chi_{\{\hat{u}> \alpha\}}.
\end{equation}

\end{thm}

In Section \ref{sec:fr} we introduce some technical machinery,
and in Section \ref{sec-mainres} prove a sequence of claims leading to the desired result. The non-locality of the operator requires new techniques in proving \eqref{cond-ineq}.




\section{A toolbox for the fractional Laplacian}\label{sec:fr}

In this section we will present a short introduction about fractional Laplace equation mainly following \cite{Silvestre_diss} and \cite{DPV}, but also some other authors cited below.

Let us first define the following fractional Sobolev spaces.  Hence, for $0<s<1$ we define
$$
H^s(\R^n)=\{v\in L^2(\R^n)\colon [v]^2_s<\infty\},
$$
where
$$
[v]^2_s = \frac{1}{2}\iint_{\R^{2n}}\frac{(v(x)-v(y))^2}{|x-y|^{n+2s}}\, dxdy.
$$
is the so-called Gagliardo-Nirenberg semi-norm.

Observe that $H^s(\R^n)$ is a Hilbert space with inner product given by
$$
(u,v)_s = \int_{\R^n} u(x)v(x)\, dx + \frac12\iint_{\R^{2n}} \frac{(u(x)-u(y))(v(x)-v(y))}{|x-y|^{n+2s}}\, dxdy.
$$
Further we define $H^{-s}(\R^n)$ as the dual space of $H^{s}(\R^n)$ and for a domain $D\subset \R^n$, 
$$
H^s_0(D)=\{v\in H^s(\R^n)\colon  v(x)=0\ \text{ in }\ D^c \}.
$$ 
Observe that $H^s_0(D)\subset H^s(\R^n)$ is a closed subspace and hence is also a Hilbert space.

We denote by $H^{-s}(D)$ the dual space of $H^s_0(D)$. Recall that if $f\in H^{-s}(\R^n)$ then, the restriction of $f$ to $H^s_0(D)$ uniquely defines a function in $H^{-s}(D)$. In that sense, we will say that $H^{-s}(\R^n)\subset H^{-s}(D)$ (even if this inclusion is not an injection).

Recall that the Gagliardo-Nirenberg semi-norm is G\^ateaux-differentiable and 
\begin{equation}\label{GN-der}
\lim_{\epsilon\to 0}\epsilon^{-1}([u+\epsilon v]^2_s-[u]^2_s) = \iint_{\R^{2n}}\frac{(u(x)-u(y))(v(x)-v(y))}{|x-y|^{n+2s}}\, dxdy.
\end{equation}


For a function $u\in H^s(\R^n)$ we can also define 
\begin{equation}\label{frac-lap}
(-\Delta)^s u(x) = p.v. \int_{\R^n}\frac{u(x)-u(y)}{|x-y|^{n+2s}}dy = \lim_{\epsilon\to 0}(-\Delta)^s_\epsilon u(x),
\end{equation}
where
$$
(-\Delta)^s_\epsilon u(x) = \int_{\R^n\setminus B_\epsilon(x)}\frac{u(x)-u(y)}{|x-y|^{n+2s}}\, dy.
$$
One can show that $(-\Delta)^s u(x)\in H^{-s}(\R^n)$, the limit in \eqref{frac-lap} holds in $H^{-s}(\R^n)$ and 
$$
\langle (-\Delta)^s u, v\rangle = \frac12\iint_{\R^{2n}} \frac{(u(x)-u(y))(v(x)-v(y))}{|x-y|^{n+2s}}\, dxdy \leq [u]_s [v]_s
$$
for any $v\in H^s(\R^n)$, where $\langle\cdot,\cdot\rangle$ is the duality product between $H^{-s}(\R^n)$ and $H^s(\R^n)$ (see \cite{DPV}).

The lemma below is the fractional analogue of the Poincar\'e inequality (see \cite[Lemma 2.4]{BLP}).
\begin{lem} \label{poincare}
Let $s \in (0,1)$, $D \in \mathbb{R}^n$ be an open and bounded set. Then we have,
\begin{equation}\label{Poincare-frac}
\| u \|^2_{L^2 (D)} \leq C(n,s,D) [u]^2_s, \ \ \ \ \text{for every} \ \ u \in H^s_0 (D),
\end{equation}
where the geometric quantity $C(n,s,D)$ is defined by
\begin{equation*}
C(n,s,D) = \min \left\{ \frac{\text{diam}\left( D \cup B  \right)^{n + 2s}}{|B|} :  B \subset \mathbb{R}^n \backslash D  \ \ \text{is a ball}  \right\}.
\end{equation*}
\end{lem}

For a function $f\in H^{-s}(D)$ we say $u_f\in H^s_0(D)$ solves the fractional boundary value problem in $D$ with homogeneous Dirichlet boundary condition
\begin{equation}\label{main-fr}
\begin{cases} (-\Delta)^s u_f(x) = f (x)  & \mbox{in } D,
 \\
 u_f=0  & \mbox{in }  D^c,  \end{cases}
\end{equation}
if 
\begin{equation}\label{partint-fr}
\frac{1}{2}\iint_{\R^{2n}}\frac{(u_f(x)-u_f(y))(v(x)-v(y))}{|x-y|^{n+2s}}dxdy=\int_D f(x)v(x)dx
\end{equation}
for any $v\in H^s_0(D)$.

The next lemma is an easy consequence of the Riesz representation Theorem, with the help of \eqref{Poincare-frac}.
\begin{lem} \label{sup}
The equation \eqref{main-fr} has a unique weak solution which satisfies
\begin{eqnarray*}
&&\int_D f u_f dx = \frac{1}{2} \iint_{\mathbb{R}^{n \times n}} \frac{\left| u_f (x) - u_f (y) \right|^2}{|x-y|^{n+2s}} dxdy\\
&=&\sup_{u \in H^s_0(D)} \left\{ 2\int_D uf - \frac{1}{2} \iint_{\mathbb{R}^{n \times n}} \frac{\left| u (x) - u (y) \right|^2}{|x-y|^{n+2s}}dxdy \right\}.
\end{eqnarray*}
\end{lem}

The following lemma can be found in \cite{Silvestre_diss}.

\begin{lem}\label{continuity}
Let $f=(-\Delta)^s u$. Assume that $f, u \in L^\infty (\mathbb{R}^n)$ and $s > 0$. Then 
\begin{enumerate}
\item If $2s \leq 1$, then $u\in C^{0,\alpha}(\mathbb{R}^n)$ for any $\alpha < 2s$. Moreover
\begin{equation*}
\| u \|_{C^{0,\alpha}(\mathbb{R}^n)} \leq  C \left(  \| u \|_{L^\infty}  + \| f \|_{L^\infty}  \right),
\end{equation*}
for a constant $C$ depending only on $n,\alpha$ and $s$.
\item If $2s > 1$, then $u\in C^{1,\alpha}(\mathbb{R}^n)$ for any $\alpha < 2s -1$. Moreover
\begin{equation*}
\| u \|_{C^{1,\alpha}(\mathbb{R}^n)} \leq  C \left(  \| u \|_{L^\infty}  + \| f \|_{L^\infty}  \right),
\end{equation*}
for a constant $C$ depending only on $n,\alpha$ and $s$.
\end{enumerate}
\end{lem}

The above results are valid also for solutions of $f=(-\Delta)^s u$ in bounded 
domains (see remarks after \cite[Proposition 2]{PSV}).

The following compactness results (see \cite[Lemma 10]{PSV}) will be used in our proofs.
\begin{lem} \label{compact}
Let $n \geq 1$, $D\in \mathbb{R}^n$ be a Lipschitz open bounded set and $\mathfrak{J}$ be a bounded subset of $L^2 (D)$. Suppose that 
\begin{equation*}
\sup_{f \in \mathfrak{J}} \int_D \int_D \frac{\left| f(x)-f(y) \right|^2}{\left| x-y  \right|^{n+2s}} dxdy < \infty.
\end{equation*}
Then, $\mathfrak{J}$ is pre-compact in $ L^2 (D)$.
\end{lem}


As a final result, we state for further reference the following lemma concerning some convex maximization problem. The proof of these facts are easy, well-known and are left to the reader.

\begin{lem}\label{maxfunc}
The set $\bar{\mathcal{R}}_\beta$ is the weak closure of the set $\mathcal{R}_\beta$. 
Moreover,
$$
\text{ext}(\bar{\mathcal{R}}_\beta)=\mathcal{R}_\beta,
$$
where $\text{ext}(C)$ denotes the extreme points of the convex set $C$.

Finally, if $g \in L^2_+ (D)$, then there exists $f \in \text{ext}(\bar{\mathcal{R}}_\beta)=\mathcal{R}_\beta$ such that
$$
\int_D hg\leq\int_D fg,
$$ 
for all $h \in \bar{\mathcal{R}}_\beta$.
\end{lem}


\section{Proof of Theorem \ref{thm:main}}\label{sec-mainres}

In this section we prove our main result, Theorem \ref{thm:main}. We will divide the proof into a series of claims.

{\bf Claim 1:} Existence.

Let
\begin{equation*}
I =  \sup_{f \in \bar{\mathcal{R}}_\beta}\int_D f u_f dx.
\end{equation*}
We first show that $I$ is finite. Consider $f \in \bar{\mathcal{R}}_\beta$. Then, by Lemma \ref{sup}  $u_f$ satisfies
\begin{equation*}
 \int_{D} f u_f\, dx = \frac{1}{2} [u_f]^2_s.
\end{equation*} 
Using H\"older's inequality and \eqref{Poincare-frac},
\begin{equation} \label{L2}
\int_D f u_f\, dx \le \| f \|_2 C [u_f]_s,
\end{equation}
and thus we obtain
\begin{equation}\label{nso}
\int_D f u_f\, dx \leq  C \| f \|^2_2 \le C,
\end{equation}
since $0\le f\le 1$ a.e. in $D$, which proves that $I$ is finite.

Let now $\{ f_i \}_{i\in\N}\subset \bar{\mathcal{R}}_\beta$ be a maximization sequence and let $u_i = u_{f_i}$. Then 
\begin{equation*}
I = \lim_{i \rightarrow \infty} \int_D f_i u_i dx.
\end{equation*}
It is clear from \eqref{L2} and \eqref{nso} that $u_i$ is bounded both in $H^s_0(D)$, hence by Lemma \ref{compact} there exist a subsequence (still denoted by $u_i$)  that converges strongly to $u_0 \in L^2 (D)$ and weakly in $H^s_0(D)$. Since $[\cdot]_s^2$ is convex, it follows that it is sequentially weakly lower semicontinuous and hence
\begin{equation} \label{ineq}
[u_0]_s^2\le \liminf_{i\to\infty}\ [u_i]_s^2 = I.
\end{equation}

On the other hand, since $f_i$ is bounded in $L^2 (D)$ and in $L^\infty(D)$, there exist  a subsequence (still denoted by $f_i$) converging weakly in $L^2(D)$ and weakly* in $L^\infty(D)$ to some $\eta \in L^\infty(D)$. Since $\bar{\mathcal{R}}_\beta$ is weakly closed, we have $\eta \in \bar{\mathcal{R}}_\beta$. Thus, we obtain 
\begin{equation} \label{remind}
\int_D f_i u_i\, dx \to \int_D \eta u_0\, dx.
\end{equation}
By Lemma \ref{sup}, \eqref{ineq} and \eqref{remind}, we obtain
\begin{eqnarray}
I &=& \lim_{i \to \infty} \int_D f_i u_i dx = \lim_{i \to \infty} 2\int_D u_i f_i\, dx - [u_i]_s^2 \\
&\leq& 2\int_D u_0\eta\, dx - [u_0]_s^2. \label{first}
\end{eqnarray}
According to Lemma \ref{maxfunc}, there exists $\hat{f} \in \mathcal{R}_\beta$ such that
\begin{equation} \label{hatoccur}
\int_D \hat{f} u_0 \, dx = \sup_{h \in \bar{\mathcal{R}}_\beta} \int_D h u_0\,  dx.
\end{equation}
Applying again Lemma \ref{sup} together with \eqref{first}, \eqref{hatoccur}, we obtain,
\begin{eqnarray*}
I & \leq &2\int_D \hat{f} u_0 \, dx - [u_0]_s^2 \\
& \leq & 2\int_D \hat{f} \hat{u} - [\hat u]_s^2\\
& = & \int_D \hat{f} \hat{u} dx \leq I,
\end{eqnarray*}
where $\hat{u} = u_{\hat{f}}$. Thus, $\hat{f}$ is a maximizer of $\Phi_s$.

\vspace{4mm}

From now on $\hat{f}  \in \bar{\mathcal{R}}_\beta$ will denote {\bf any} maximizer of $\Phi_s$, not necessary the one 
obtained in Claim 1, which we already know belongs to $ \mathcal{R}_\beta$.

\vspace{4mm}

{\bf Claim 2}:  $\hat{f}$ maximizes the linear functional $L(f):= \int_D \hat{u} f dx$ over $\bar{\mathcal{R}}_\beta$.
 
\vspace{2mm}

Let us take $ f\in \bar{\mathcal{R}}_\beta$ and use the maximization property
$$
 \Phi_s((1-\epsilon)\hat{f}+\epsilon f)\leq\Phi_s(\hat{f}).
$$
This inequality implies that
$$
\epsilon \iint_{\R^{2n}} \frac{(\hat u(x)-\hat u(y))((u_f-\hat u)(x) - (u_f-\hat u)(y))}{|x-y|^{n+2s}}\, dxdy + 2\epsilon^2 [u_f - \hat u]_s^2\le 0.
$$
If we now divide by $\epsilon$ and take the limit as $\epsilon\to 0$ we get
$$
\frac12 \iint_{\R^{2n}} \frac{(u_f(x)-u_f(y))(\hat u(x)-\hat u(y))}{|x-y|^{n+2s}}\, dxdy \le \frac12[\hat u]_s^2.
$$

But if we now use Lemma \ref{sup}, this last inequality becomes
$$
\int_D f\hat{u}\, dx \leq \int_D \hat{f}\hat{u}\, dx,
$$
as we wanted to show.

\vspace{4mm}


Next, observe that from  Lemma \ref{maxfunc}, there exists a $\tilde{f} = \chi_{E} \in \mathcal{R}_\beta=\text{ext}(\bar{\mathcal{R}}_\beta)$ such that $\tilde{f}$ maximizes $L(f)$ over $\bar{\mathcal{R}}_\beta$.
\vspace{4mm}

{\bf Claim 3}: $\alpha = \sup_{x \in E^c} \hat{u} (x)\leq \gamma = \inf_{x \in E} \hat{u} (x)$ (where $\sup$ and $\inf$ denote the essential supremum and the essential infimum respectively).

\vspace{2mm}

Assume by contradiction that $\gamma < \alpha$. Let us fix $\gamma < \xi_1 < \xi_2 < \alpha.$ Since $\xi_1 > \gamma$, there exists a set $A \in E$, with positive measure, such that $\hat{u} \leq \xi_1$ on $E$. Similarly, $\xi_2 < \alpha$ implies that there exists a $B \in E^c$, with positive measure, such that $\hat{u} \geq \xi_2$ on $E^c$. Without loss of generality, we assume that $A$ and $B$ have the same Lebesgue measure. Next,  we define a new rearrangement of $\tilde{f}$, which is denoted by $\bar{f}\in  \mathcal{R}_\beta$.
\begin{equation*}
\bar{f} = 
\begin{cases}
0, & x\in A;\\
1, & x\in B;\\
\tilde{f}(x), & x\in D \backslash \left( A \cup B \right).
\end{cases}
\end{equation*}
Therefore,
\begin{eqnarray*}
\int_D \bar{f} \hat{u} \, dx - \int_D \tilde{f} \hat{u} \, dx &=&  \int_B \bar{f} \hat{u} \, dx - \int_A \tilde{f} \hat{u}\, dx \\
&\geq & \xi_2 \int_B \bar{f} \, dx - \xi_1 \int_A \tilde{f} \, dx \\
&=& \left( \xi_2 - \xi_1 \right) \int_A \tilde{f} \, dx > 0,
\end{eqnarray*}
which contradicts the maximality of $\hat{f}$.

Recall that $\hat u$ is continuous (Lemma \ref{continuity}), therefore $\alpha=\gamma$.

\vspace{4mm}

{\bf Claim 4}: $\chi_{\{\hat u>\alpha\}}\le \hat f\le \chi_{\{\hat u\ge\alpha\}}$.

\vspace{2mm}

We need to prove that
\begin{equation*}
\hat{f}=
\begin{cases}
1 \ \ a.e. \ \text{in} \ \left\{ \hat{u} > \alpha  \right\};\\
0 \ \ a.e. \ \text{in} \ \left\{ \hat{u} < \alpha  \right\}.
\end{cases}
\end{equation*}

We argue by contradiction. Assume there exists a $A \subset \left\{ \hat{u} > \alpha  \right\}$, with positive measure, such that $\hat{f} < 1$ in $A$. Since $|\{ \hat{u} > \alpha \}| \leq \beta$, $\hat{f} > 0$ in some subset of $\{\hat{u} \leq \alpha  \}$. Thus, we can replace the function $\hat{f}$ by a function $f \in \bar{\mathcal{R}}_\beta$ which has larger values in $A$ and smaller values in $\{\hat{u} \leq \alpha  \}$. As a result,
\begin{equation*}
\int_D f \hat{u} dx > \int_D \hat{f} \hat{u} dx,
\end{equation*}
which contradicts the maximality of $\hat{f}$. Therefore, $\hat{f} = 1$ a.e. in $\{ \hat{u} > \alpha \}$.

Similarly, assume there exists a $A \subset \{ \hat{u} < \alpha \}$, with positive measure such that $\hat{f} > 0$ in $A$. Since $E \subset \{ \hat{u} \geq \alpha \}$, $\hat{f} < 1$ in some subset of $\{ \hat{u} \geq \alpha \}$. Thus, we can replace the function $\hat{f}$ by a function $f \in \bar{\mathcal{R}}_\beta$ which vanishes in $A$ and has larger values in $\{ \hat{u} \geq \alpha \}$. As a result,
\begin{equation*}
\int_D f \hat{u} dx > \int_D \hat{f} \hat{u} dx,
\end{equation*}
which contradicts the maximality of $\hat{f}$. Therefore, $\hat{f} = 0$ a.e. in $\{ \hat{u} < \alpha \}$.

\vspace{4mm}

{\bf Claim 5}: $|\{\hat{u} = \alpha\}| = 0.$

\vspace{2mm}

Assume $|\{\hat{u} = \alpha\}| > 0$. Take $\tilde{E} \subset D$ such that $\{ \hat{u} > \alpha \} \subset \tilde{E} \subset \{ \hat{u} \geq \alpha \}$ and $|\tilde{E}| = \beta$. Let $v \in H^s_0 (D)$ be the unique solution to the following fractional boundary value problem,
\begin{equation*}
\begin{cases}
(-\Delta^s v) = \chi_{\tilde{E}} &  \text{in } D,\\
v = 0  &   \text{in }   D^c.
\end{cases}
\end{equation*}
Set $\tilde{u} := \frac12\hat{u} + \frac12 v$. Then $(-\Delta)^s \tilde{u} = \frac12\hat{f} + \frac12 \chi_{\tilde{E}} \in \bar{\mathcal{R}}_\beta$. Now, it suffices to show that 
\begin{equation}\label{biggerzero}
\left[ \tilde{u} \right]^2_s  > \left[ \hat{u}  \right]^2_s,
\end{equation}
which would contradict the maximality of $\hat{u}$. But, by elementary computations, \eqref{biggerzero} is equivalent to
\begin{equation}\label{keycontr}
\frac12 \left[ \hat{u} - v \right]^2_s > 2 \iint_{\R^{2n}} \frac{(\hat u(x)-\hat u(y))((\hat u - v)(x) - (\hat u - v)(y))}{|x-y|^{n+2s}}\, dxdy.
\end{equation}

Next, from Lemma \ref{sup} and Claim 4, we get
\begin{align*}
\iint_{\R^{2n}} \frac{(\hat u(x)-\hat u(y))((\hat u - v)(x) - (\hat u - v)(y))}{|x-y|^{n+2s}}&\, dxdy\\
& =  2 \int_D \hat{u} \left( \hat{f} - \chi_{\tilde{E}} \right) \, dx \\
&= 2\alpha \int_{\{\hat{u}=\alpha\}} \left(\hat{f} - \chi_{\tilde{E}}\right)\, dx\\
&= 2\alpha \int_{D} \left(\hat{f} - \chi_{\tilde{E}}\right)\, dx\\
&= 2\alpha (\beta - \beta) =0.
\end{align*}
This completes the proof of the claim.

\vspace{4mm}

Claims 4 and 5 imply
$$ 
\hat{f}=\chi_{\{\hat{u}>\alpha  \}}
$$
and the proof of Theorem \ref{thm:main} is complete.

\vspace{4mm}

\begin{rem}
As in the classical case, it is in general {\bf not} true that the function $\hat{u}(x) $ minimizes the (non-convex) functional
\begin{equation}\label{func-fr}
J(u)=[u]_s^2 -2\int_D \chi_{\{u>\alpha\}}u  dx,
\end{equation}
over $ H_0^s(D)$. 
  \end{rem}

\begin{proof}
Let us first introduce the subset of functions which do not have flat positive components as follows
$$
\tilde{H}_0^s(D)=\{u\in H_0^s(D)\,\,|\,\,\mathcal{L}_N (\hat{u} = t) = 0 \text{ for all } t>0\}.
$$
Since $\tilde{H}_0^s(D)$ is dense in $ H_0^s(D)$ we can replace $H_0^s(D)$ by $\tilde{H}_0^s(D)$ while taking supremum or infimum. Using the fact that for a function $u\in \tilde{H}_0^s(D)$
we can always find a real number $\alpha_u$ such that $|\{u>\alpha_u\}|=\beta$, we obtain
\begin{multline}
\Phi_s(\hat{f})=\max_{f\in \bar{\mathcal{R}}_\beta}\sup_{u\in H_0^s(D)} \left(2\int_D fudx -[u]_s^2\right)=\\
\sup_{u\in \tilde{H}_0^s(D)}\sup_{f\in \bar{\mathcal{R}}_\beta} \left(2\int_D fudx -[u]_s^2\right)=
\sup_{u\in \tilde{H}_0^s(D)} \left(2\int_D \chi_{\{u>\alpha_u\}}udx -[u]_s^2\right)=
\\
-\inf_{u\in \tilde{H}_0^s(D)} \left([u]_s^2-2\int_D \chi_{\{u>\alpha_u\}}udx \right),
\end{multline}
which implies that 
$$
[\hat{u}]_s^2-2\int_D\hat{f}\hat{u}dx=
[\hat{u}]_s^2-2\int_D
\chi_{\{{\hat u}>\alpha\}}
\hat{u}dx=\inf_{u\in \tilde{H}_0^s(D)} \left([u]_s^2-2\int_D \chi_{\{u>\alpha_u\}}udx \right).
$$
However 
$$
[\hat{u}]_s^2-2\int_D
\chi_{\{{\hat u}>\alpha\}}
\hat{u}dx\not= \inf_{u\in H_0^s(D)} \left([u]_s^2-2\int_D \chi_{\{u>\alpha\}}udx \right).
$$
A simple heuristic example can be observed as follows. Consider $D $ which consists of two disconnected balls.
We can always connect them by a very narrow tube, which would preserve the discussion below unchanged. 
For small values of $\beta$
the maximizer of the optimal rearrangement problem will concentrate the set $\{\hat{u}>\alpha\}$ in one of the two balls and keep the function zero in the other ball. On contrast 
the minimizer of the right hand side can reach a smaller value by ``copying'' the non-zero function to the ball where $\hat{u}$
is zero. 
\end{proof}

\subsection*{Acknowledgment}
The research of Zhiwei Cheng and Hayk Mikayelyan has been partly supported by the National Science Foundation of China (grant no.1161101064). Juli\'an F. Bonder is supported by by grants UBACyT UBACYT Prog. 2018 20020170100445BA, CONICET PIP 11220150100032CO and ANPCyT PICT 2016-1022. 
 
\bibliographystyle{plain}
\bibliography{rearr}
\end{document}